\documentclass{article}

\usepackage{amssymb}

\usepackage{graphicx}

\title{Some Results
on Lagrangians of  Hypergraphs
}

\author{ Qingsong Tang \thanks{Mathematics School, Institute of Jilin University, Changchun, 130012, China, and College of Sciences, Northeastern University, Shenyang, 110819, China.Email: t\_qsong@sina.com.cn}  \and Yuejian Peng \thanks{ School of Mathematics, Hunan University, Changsha 410082, P.R. China and Indiana State University, Terre Haute, IN, 47809, USA. Email: ypeng1@163.com}   \and  Xiangde Zhang \thanks{College of Sciences, Northeastern University, Shenyang, 110819, China} \and Cheng Zhao \thanks{Department of Mathematics and Computer Science, Indiana State University, Terre Haute, IN, 47809 and School of Mathematics, Jilin University, Changchun 130012, P.R. China. Email: cheng.zhao@indstate.edu}}

\date{}

\newtheorem{defi}{Definition}[section]
\newtheorem{theo}{Theorem}[section]

\newtheorem{remark}[theo]{Remark}

\newtheorem{lemma}[theo]{Lemma}

\newtheorem{coro}[theo]{Corollary}
\newtheorem{con}[theo]{Conjecture}

\newtheorem{fact}[theo]{Fact}

\newcommand{\qed}{\hspace*{\fill} \rule{7pt}{7pt}}

\begin{document}
\maketitle
\begin{abstract}
In 1965, Motzkin and Straus \cite{MS} provided a new proof of Tur\'an's theorem based on a continuous characterization of the clique number of a graph using the Lagrangian of a graph. This new proof aroused interests in the study of Lagrangians of $r$-uniform graphs. The Lagrangian of a hypergraph has been  a useful tool in hypergraph extremal problems. Sidorenko  and Frankl-F\"{u}redi  applied Lagrangians of hypergraphs in finding Tur\'an densities of hypergraphs.
Frankl and R\"odl  applied it in disproving Erd\"os'  jumping constant conjecture. In most applications, we need an upper bound for the Lagrangian of a hypergraph.
Frankl and F\"{u}redi conjectured that the $r$-uniform graph with $m$ edges formed by taking the first $m$ sets in the colex ordering of ${\mathbb N}^{(r)}$ has the largest Lagrangian of all $r$-uniform graphs with $m$ edges. Talbot in \cite{T} provided some evidence for Frankl and F\"{u}redi's conjecture.  In this paper, we prove that  the $r$-uniform graph with $m$ edges formed by taking the first $m$ sets in the colex ordering of ${\mathbb N}^{(r)}$ has the largest Lagrangian of all $r$-uniform graphs on $t$ vertices with $m$ edges when $m={t\choose r}-3$ or ${t\choose r}-4$. As an implication, we also prove that Frankl and F\"{u}redi's conjecture holds for 3-uniform graphs with $m={t\choose 3}-3$ or ${t\choose 3}-4$ edges.
\end{abstract}

Key Words:  Lagrangians of $r$-uniform graphs, Extremal problems in hypergraphs

\section{Introduction}

In 1941, Tur\'an \cite{Turan} provided an answer to the following question: What is the maximum number of edges in a graph on $n$ vertices without containing a complete subgraph of order $t$, for a given $t$? This is the well-known Tur\'an theorem. Later, in another classical paper, Motzkin and Straus \cite{MS} provided a new proof of Tur\'an's theorem based on a continuous characterization of the clique number of a graph using Lagrangians of  graphs. This new proof aroused interests in the study of Lagrangians of $r$-graphs. The Lagrangian of a hypergraph has been  a useful tool in hypergraph extremal problems. For example, Sidorenko \cite{sidorenko89} and Frankl-F\"{u}redi \cite{FF} applied Lagrangians of hypergraphs in finding Tur\'an densities of hypergraphs.
Frankl and R\"odl \cite{FR84} applied it in disproving Erd\"os' jumping constant conjecture.  More applications of Lagrangians can be found in \cite{sidorenko93}, \cite{FR89} and \cite{mubayi06}. In most applications, we need an upper bound for the Lagrangian of a hypergraph. In the course of estimating Tur\'an densities of some hypergraphs, Frankl and F\"uredi \cite{FF} asked the following question: Given $r \ge 3$ and $m \in{\mathbb N}$, how large can the Lagrangian of an $r$-graph with $m$ edges be? Before stating their conjecture on this problem, we give some definitions and notation.

For a set $V$ and a positive integer $r$ we denote by $V^{(r)}$ the family of all $r$-subsets of $V$. An $r$-uniform graph or $r$-graph $G$ consists of a set $V(G)$ of vertices and a set $E(G) \subseteq V(G) ^{(r)}$ of edges. An edge $e=\{a_1, a_2, \ldots, a_r\}$ will be simply denoted by $a_1a_2 \ldots a_r$. An $r$-graph $H$ is  a {\it subgraph} of an $r$-graph $G$, denoted by $H\subseteq G$ if $V(H)\subseteq V(G)$ and $E(H)\subseteq E(G)$.  Let ${\mathbb N}$ be the set of all positive integers. For any  $n \in {\mathbb N}$, we denote the set $\{1, 2, 3, \ldots, n\}$ by $[n]$. Let $K^{(r)}_t$ denote the complete $r$-graph on $t$ vertices, that is the $r$-graph on $t$ vertices containing all possible edges. A complete $r$-graph on $t$ vertices is also called a clique with order $t$. We also let $[n]^{(r)}$  represent the  complete $r$-uniform graph on the vertex set $[n]$. When $r=2$, an $r$-uniform graph is a simple graph.  When $r\ge 3$,  an $r$-graph is often called a hypergraph.

For an $r$-graph $G=(V,E)$ we denote the $(r-1)$-neighborhood of a vertex $i \in V$ by $E_i=\{A \in V^{(r-1)}: A \cup \{i\} \in E\}$. Similarly, we will denote the $(r-2)$-neighborhood of a pair of vertices $i,j \in V$ by $E_{ij}=\{B \in V^{(r-2)}: B \cup \{i,j\} \in E\}$. We denote the complement of $E_i$ by $E^c_i=\{A \in V^{(r-1)}: A \cup \{i\} \in V^{(r)} \backslash E\}$. Also, we  denote the complement of $E_{ij}$ by
$E^c_{ij}=\{B \in V^{(r-2)}: B \cup \{i,j\} \in V^{(r)} \backslash E\}$. Denote $$E_{i\setminus j}=E_i\cap E^c_j.$$


\begin{defi}
For  an $r$-uniform graph $G$ with the vertex set $[n]$,
edge set $E(G)$ and a vector $\vec{x}=(x_1,\ldots,x_n) \in R^n$,
define
$$\lambda (G,\vec{x})=\sum_{i_1i_2 \cdots i_r \in E(G)}x_{i_1}x_{i_2}\ldots x_{i_r}.$$

Let $S=\{\vec{x}=(x_1,x_2,\ldots ,x_n): \sum_{i=1}^{n} x_i =1, x_i
\ge 0 {\rm \ for \ } i=1,2,\ldots , n \}$. The Lagrangian of
$G$, denote by $\lambda (G)$, is defined as
 $$\lambda (G) = \max \{\lambda (G, \vec{x}): \vec{x} \in S \}.$$

We call $\vec{x}=(x_1, x_2, \ldots, x_n) \in R^n$ a legal weighting for $G$ if
$\vec{x}\in S$. A vector $\vec{y}\in S$ is called an {\em optimal weighting} for $G$ if $\lambda (G, \vec{y})=\lambda(G)$.
\end{defi}
 The following fact is easily implied by the definition of the Lagrangian.

\begin{fact}\label{mono}
Let $G_1$, $G_2$ be $r$-uniform graphs and $G_1\subseteq G_2$. Then $\lambda (G_1) \le \lambda (G_2).$
\end{fact}

In \cite{MS}, Motzkin and Straus provided the following simple expression for the Lagrangian of a 2-graph.

\begin{theo} (Motzkin and Straus \cite{MS}) \label{MStheo}
If $G$ is a 2-graph in which a largest clique has order $t$ then
$\lambda(G)=\lambda(K^{(2)}_t)={1 \over 2}(1 - {1 \over t})$.
\end{theo}

An attempt to generalize the Motzkin-Straus theorem to hypergraphs is due to S\'os and Straus \cite{SS}.  Recently, in \cite{BP1} and \cite{BP2} Rota Bul\'o and Pelillo generalized the Motzkin and Straus' result to $r$-graphs in some way using a continuous characterization of maximal cliques.
Determining the Lagrangian of a general $r$-graph is non-trivial when $r \ge 3$. Indeed the obvious generalization of Motzkin and Straus' result is false because there are many examples of $r$-graphs that do not achieve their Lagrangian on any proper subhypergraph.

For distinct $A, B \in {\mathbb N}^{(r)}$ we say that $A$ is less than $B$ in the {\em colex ordering} if $max(A \triangle B) \in B$, where $A \triangle B=(A \setminus B)\cup (B \setminus A)$. For example we have $246 < 156$ in ${\mathbb N}^{(3)}$ since $max(\{2,4,6\} \triangle \{1,5,6\}) \in \{1,5,6\}$. In colex ordering, $123<124<134<234<125<135<235<145<245<345<126<136<236<146<246<346<156<256<356<456<127<\cdots .$ The following conjecture of Frankl and F\"uredi (if it is true) proposes a  solution to the question mentioned at the beginning.

\begin{con} (Frankl and F\"uredi \cite{FF})\label{conjecture} The $r$-graph with $m$ edges formed by taking the first $m$ sets in the colex ordering of ${\mathbb N}^{(r)}$ has the largest Lagrangian of all $r$-graphs with  $m$ edges. In particular, the $r$-graph with $t \choose r$ edges and the largest Lagrangian is $[t]^{(r)}$.
\end{con}

This conjecture is true when $r=2$ by Theorem \ref{MStheo}. For the case $r=3$, Talbot in \cite{T} proved the following.

\begin{theo} (Talbot \cite{T}) \label{Tal} Let $m$ and $t$ be integers satisfying
$${t-1 \choose 3} \le m \le {t-1 \choose 3} + {t-2 \choose 2} - (t-1).$$
Then Conjecture \ref{conjecture} is true for $r=3$ and this value of $m$.  Conjecture \ref{conjecture} is also true for $r=3$ and $m= {t \choose 3}-1$ or $m={t \choose 3} -2$.
\end{theo}

The truth of Frankl and F\"uredi's conjecture is not known in general for $r \ge 4$. Even in the case $r=3$, Theorem \ref{Tal} does not cover the case when ${t-1 \choose 3}+{t-2 \choose 2}-(t-2) \le m \le {t \choose 3}-3$ in this conjecture. In \cite{HPZ}, He, Peng, and Zhao verified Frankl and F\"uredi's conjecture  for small values $m$ when $r=3$. In \cite{PZ}, Peng and Zhao generalized Theorem \ref{Tal} further. Talbot in \cite{T} proved some result as evidence of the truth of Conjecture \ref{conjecture} for $r$-graphs supported on $t+1$ vertices with $m={t \choose r}$ edges. Also, in \cite{PTZ} Peng, Tang, and Zhao provided more results for Conjecture \ref{conjecture} when a given $r$-graph on $t$ vertices satisfies some conditions.

Let $C_{r,m}$ denote the $r$-graph with $m$ edges formed by taking the first $m$ sets in the colex ordering of ${\mathbb N}^{(r)}$.

\begin{lemma}  \cite{T} \label{LemmaTal7}
For  integers $m,  t, $ and $r$ satisfying ${t-1 \choose r} \le m \le {t-1 \choose r} + {t-2 \choose r-1}$,
we have $\lambda(C_{r,m}) = \lambda([t-1]^{(r)})$.
\end{lemma}

\begin{defi}
An $r$-graph $G=(V,E)$  on the vertex set $[n]$ is {\it left-compressed} if $j_1j_2 \cdots j_r \in E$ implies $i_1i_2 \cdots i_r \in E$ provided $i_p \le j_p$ for every $p, 1\le p\le r$. Equivalently,  $G$ is {\it left-compressed} if $E_{j\setminus i}=\emptyset$ for any $1\le i<j\le n$.
\end{defi}


Denote $$\lambda_{m}^{r}=\max\{\lambda(G): G {\rm \ is \ an \ } r-{\rm graph\ with \ } m {\rm \ edges }\}.$$ The following lemma implies that we only need to consider left-compressed $r$-graphs when Conjecture \ref{conjecture} is discussed.
\begin{lemma}  \cite{T} \label{LemmaTal8}
Let $m,t$ be positive integers satisfying  $m\leq{t \choose r}-1$, then there exists a left-compressed $r$-graph $G$ with $m$ edges such that $\lambda(G)=\lambda_{m}^{r}.$
\end{lemma}
In this paper, we show that
\begin{theo} \label{theorem 2} Let $m$ and $t$ be positive integers
 satisfying $m={t \choose r}-3$. Then the $r$-graph with $m$ edges formed by taking the first $m$ sets in the colex ordering of ${\mathbb N}^{(r)}$ has the largest Lagrangian of all $r$-graphs  with $m$ edges and $t$ vertices.
\end{theo}
Combining with a result in \cite{T} (Lemma \ref{LemmaTal9}), we show   the following corollary for $r=3$.
\begin{coro} \label{corollary 1} Let $m$ and $t$ be positive integers
 satisfying $m={t \choose 3}-3$. Then Conjecture \ref{conjecture} is true for $r=3$ and this value of $m$.
\end{coro}

\begin{theo} \label{theorem 3} Let $m$ and $t$ be positive integers
 satisfying $m={t \choose r}-4$. Then the $r$-graph with $m$ edges formed by taking the first $m$ sets in the colex ordering of ${\mathbb N}^{(r)}$ has the largest Lagrangian of all $r$-graphs  with   $m$ edges and $t$ vertices.
\end{theo}
Combining with a result in \cite{T} (Lemma \ref{LemmaTal9}), we show   the following corollary for $r=3$.
\begin{coro} \label{corollary 2} Let $m$ and $t$ be positive integers
 satisfying $m={t \choose 3}-4$. Then Conjecture \ref{conjecture} is true for $r=3$ and this value of $m$.
\end{coro}

The proof of Theorem \ref{theorem 2} and Corollary \ref{corollary 1} will be given in Section \ref{proof1} and the proof of  Theorem \ref{theorem 3} and  Corollary \ref{corollary 2} will be given in Section \ref{proof2}. Next, we state some useful results.

\section{Useful Results}

We will impose one additional condition on any optimal weighting ${\vec x}=(x_1, x_2, \ldots, x_n)$ for an $r$-graph $G$:
\begin{eqnarray}
 &&|\{i : x_i > 0 \}|{\rm \ is \ minimal, i.e. \ if}  \ \vec y {\rm \ is \ a \ legal \ weighting \ for \ } G  {\rm \ satisfying }\nonumber \\
 &&|\{i : y_i > 0 \}| < |\{i : x_i > 0 \}|,  {\rm \  then \ } \lambda (G, {\vec y}) < \lambda(G) \label{conditionb}.
\end{eqnarray}

When the theory of Lagrange multipliers is applied to find the optimum of $\lambda(G)$, subject to $\sum_{i=1}^n x_i =1$, notice that $\lambda (E_i, {\vec x})$ corresponds to the partial derivative of  $\lambda(G, \vec x)$ with respect to $x_i$. The following lemma gives some necessary condition of an optimal weighting of  $G$.

\begin{lemma} (Frankl and R\"odl \cite{FR84}) \label{LemmaTal5} Let $G=(V,E)$ be an $r$-graph on the vertex set $[n]$ and ${\vec x}=(x_1, x_2, \ldots, x_n)$ be an optimal  weighting for $G$ with $k$  ($\le n$) non-zero weights $x_1$, $x_2$, $\cdots$, $x_k$ satisfying condition (\ref{conditionb}). Then for every $\{i, j\} \in [k]^{(2)}$, (a) $\lambda (E_i, {\vec x})=\lambda (E_j, \vec{x})=r\lambda(G)$, (b) there is an edge in $E$ containing both $i$ and $j$.
\end{lemma}

\begin{remark}\label{r1} (a) In Lemma \ref{LemmaTal5}, part(a) implies that
$$x_j\lambda(E_{ij}, {\vec x})+\lambda (E_{i\setminus j}, {\vec x})=x_i\lambda(E_{ij}, {\vec x})+\lambda (E_{j\setminus i}, {\vec x}).$$
In particular, if $G$ is left-compressed, then
$$(x_i-x_j)\lambda(E_{ij}, {\vec x})=\lambda (E_{i\setminus j}, {\vec x})$$
for any $i, j$ satisfying $1\le i<j\le k$ since $E_{j\setminus i}=\emptyset$.

(b) If  $G$ is left-compressed, then for any $i, j$ satisfying $1\le i<j\le k$,
\begin{equation}\label{enbhd}
x_i-x_j={\lambda (E_{i\setminus j}, {\vec x}) \over \lambda(E_{ij}, {\vec x})}
\end{equation}
holds.  If  $G$ is left-compressed and  $E_{i\setminus j}=\emptyset$ for $i, j$ satisfying $1\le i<j\le k$, then $x_i=x_j$.

(c) By (\ref{enbhd}), if  $G$ is left-compressed, then an optimal weighting  ${\vec x}=(x_1, x_2, \ldots, x_n)$ for $G$  must satisfy
\begin{equation}\label{conditiona}
x_1 \ge x_2 \ge \ldots \ge x_n \ge 0.
\end{equation}
\end{remark}

\section{Proofs of Theorem \ref{theorem 2} and  Corollary \ref{corollary 1}}\label{proof1}

Denote
$$\lambda_{(m,t)}^{r}=\max\{\lambda(G): G {\rm \ is \ an \ } r{\rm -graph  \ with \ } m {\rm \ edges \ and \ } t \ {\rm vertices} \}.$$

An $r$-tuple  $i_1 i_2\cdots i_r$ is called a {\it descendant } of an $r$-tuple  $j_1j_2\cdots j_r$ if $i_s\le j_s$ for each $1\le s\le r$, and $i_1+i_2+\cdots +i_r < j_1+j_2+\cdots +j_r$. In this case, the $r$-tuple $j_1j_2\cdots j_r$   is called an {\it ancestor} of $i_1 i_2\cdots i_r$.  The $r$-tuple $i_1i_2\cdots i_r$   is called a {\it direct descendant} of $j_1 j_2\cdots j_r$ if   $i_1i_2\cdots i_r$   is a  descendant of $j_1j_2\cdots j_r$ and $j_1+j_2+\cdots +j_r=i_1+i_2+\cdots +i_r +1$.  We say that $i_1 i_2\cdots i_r$ has lower hierarchy than  $j_1j_2\cdots j_r$ if $i_1 i_2\cdots i_r$  is  a descendant of $j_1j_2\cdots j_r$. This is a partial order on the set of all $r$-tuples.  Figure 1 is a Hessian diagram on all $r$-tuples on $[t]$. In this diagram, $i_1 i_2\cdots i_r$ and $j_1j_2\cdots j_r$ are connected by an edge if $i_1i_2\cdots i_r$   is  a  direct descendant of $j_1j_2\cdots j_r$.


\begin{remark}\label{releftcom}
An $r$-graph $G$ is left-compressed if and only if all descendants of an edge of $G$ are edges of $G$. Equivalently, if an $r$-tuple is not an edge of $G$, then none of its ancestors will be an edge of $G$.
\end{remark}

\begin{lemma}\label{leftcom} There exists a left-compressed $r$-graph $G$ on the vertex set $[t]$ with $m$ edges such that $\lambda(G)=\lambda_{(m,t)}^{r}$.
\end{lemma}

{\em Proof.} Let $G'=(V, E)$ be an $r$-graph on the vertex set $[t]$ with $m$ edges  such that $\lambda(G')=\lambda_{(m,t)}^{r}$. We call such an $r$-graph $G'$ an extremal $r$-graph for $m$ and $t$. Let ${\vec x}=(x_1, x_2, \ldots, x_t)$ be an optimal weighting of $G'$. We can assume that $x_i\ge x_j$ when $i<j$ since otherwise we can just relabel the vertices of $G'$ and obtain another  extremal $r$-graph for $m$ and $t$ with an optimal weighting ${\vec x}=(x_1, x_2, \ldots, x_t)$ satisfying $x_i\ge x_j$ when $i<j$. If $G'$ is not left-compressed, then there is an edge such that at least one of its descendants is not an edge. Replace all those edges by its available descendants with the lowest hierarchy, then we get a left-compressed r-graph $G$ on  the vertex set $[t]$ with $m$ edges and $\lambda(G, {\vec x})\ge \lambda(G')$. Therefore, $G$ is a left-compressed extremal $r$-graph
for $m$ and $t$.
\qed

\bigskip
\begin{figure}
\centering
\includegraphics{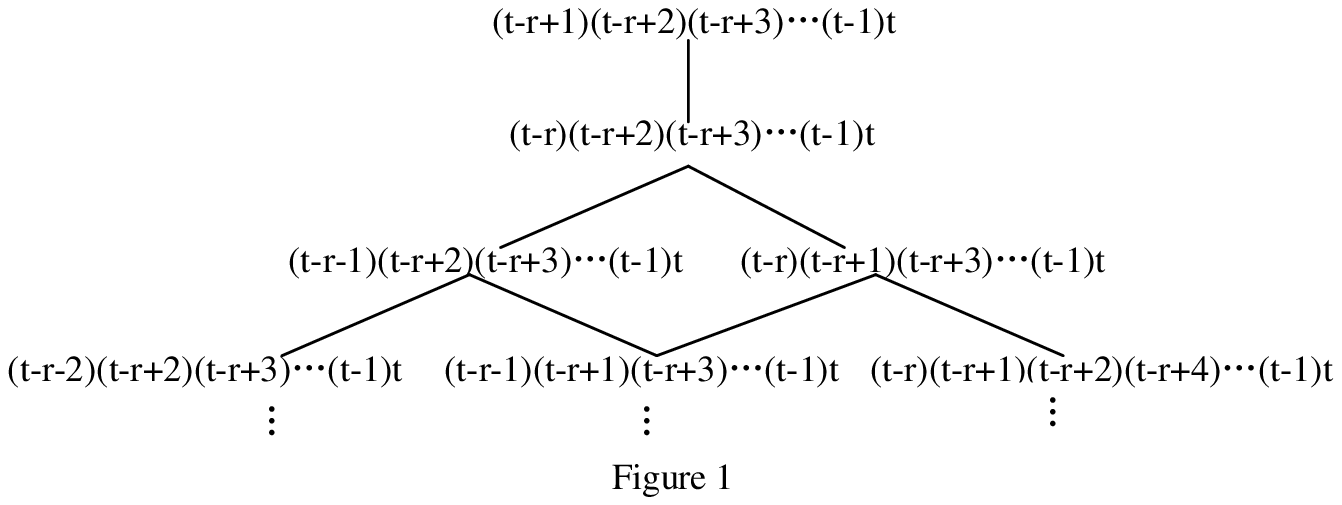}
\end{figure}
{\em Proof of Theorem \ref{theorem 2}.} By Lemma \ref{leftcom},
we only need to consider left-compressed $r$-graphs on vertex set $[t]$  with $m={t \choose r}-3$ edges. Every left-compressed $r$-graph on $[t]$  with $m={t \choose r}-3$ edges can be obtained by removing three $r$-tuples from $[t]^{(r)}$ such that if an $r$-tuple is removed then all its ancestors should be removed by Remark \ref{releftcom}. In view of Figure 1, these three $r$-tuples to be removed are either $$\{(t-r+1)(t-r+2)\ldots (t-1)t,(t-r)(t-r+2)\ldots (t-1)t,(t-r)(t-r+1)(t-r+3)\ldots (t-1)t\}$$ or $$\{(t-r+1)(t-r+2)\ldots (t-1)t,(t-r)(t-r+2)\ldots (t-1)t,(t-r-1)(t-r+2)\ldots (t-1)t\}.$$
Therefore,  there are only two different left-compressed $r$-graphs with $m={t \choose r}-3$ edges on $[t]$. They are
$$G_{1}=([t],E) \  {\rm with \ the \ edge \ set } $$
$$E=[t]^{(r)}\setminus\{(t-r+1)(t-r+2)\ldots (t-1)t,(t-r)(t-r+2)\ldots (t-1)t,(t-r)(t-r+1)(t-r+3)\ldots (t-1)t\},$$
and
 $$G_{2}=([t],E') \ {\rm with \ the \ edge \ set } $$
$$E'=[t]^{(r)}\setminus\{(t-r+1)(t-r+2)\ldots (t-1)t,(t-r)(t-r+2)\ldots (t-1)t,(t-r-1)(t-r+2)\ldots (t-1)t\}.$$
Clearly, $G_{2}$ is formed by taking the first $m$ sets in the colex ordering of ${\mathbb N}^{(r)}$. So in order to prove Theorem \ref{theorem 2}, we only need to prove $\lambda(G_{1})\leq \lambda(G_{2})$.

 Let $\vec{x}=(x_{1},x_{2},\ldots ,x_{t})$  be an optimal weighting for $G_{1}$ satisfying $x_1 \ge x_2 \ge \ldots \ge x_t \ge 0$. First note that  $x_t>0$.  If $x_t=0$, then $\lambda(G_1)=\lambda([t-1]^{(r)})$. However, if we take
a legal weighting $\vec{x}=(x_1,\ldots,x_t)$, where $x_1=x_2=\cdots=x_{t-2}={1 \over t-1}$ and $x_{t-1}=x_{t}={1 \over 2(t-1)}$, then $\lambda(G_1,\vec{x} )> \lambda([t-1]^{(r)})$. This contradiction implies that $x_t>0$.
Since $G_{1}$ is left-compressed and $E_{i\setminus j}=\emptyset$ for $i, j$ satisfying $1\le i<j\le t-r-1$, or $t-r\le i<j\le t-r+2$ or $t-r+3\le i<j\le t$,
by Remark \ref{r1}(b), we have $x_1=x_2=\cdots=x_{t-r-1}=a$, $x_{t-r}=x_{t-r+1}=x_{t-r+2}=b$, and $x_{t-r+3}=x_{t-r+4}=\cdots=x_{t}=c$.

Consider a  weighting for $G_{2}$:  ${\vec y}=(y_1, y_2, \ldots, y_t)$ given by $y_i=x_i$ for $i\neq t-r-1$, $i\neq t-r$ and $y_{t-r-1}=x_{t-r-1}-\delta$, $y_{t-r}=x_{t-r}+\delta$. Then
\begin{eqnarray}\label{eq10}
\lambda(G_{2},\vec{y})-\lambda(G_{2},\vec{x})&=& \delta[\lambda(E'_{t-r},\vec{x})-\lambda(E'_{t-r-1},\vec{x})]-\delta^{2}\lambda(E'_{(t-r-1)(t-r)},\vec{x}) \nonumber\\
&=&\delta(x_{t-r-1}-x_{t-r})\lambda(E'_{(t-r-1)(t-r)},\vec{x})-\delta^{2}\lambda(E'_{(t-r-1)(t-r)},\vec{x})\nonumber \\
&=&\delta(a-b)\lambda(E'_{(t-r-1)(t-r)},\vec{x})-\delta^{2}\lambda(E'_{(t-r-1)(t-r)},\vec{x}).
\end{eqnarray}
Let $\delta=\frac{a-b}{2}$, then $y_{t-r-1}=a-\frac{a-b}{2}=\frac{a+b}{2}>0$, and $y_{t-r}=b+\frac{a-b}{2}=\frac{a+b}{2}$. Hence ${\vec y}=(y_1, y_2, \ldots, y_t)$ is also a legal weighting for $G_{2}$, and
 $$\lambda(G_{2},\vec{y})-\lambda(G_{2},\vec{x})=\frac{(a-b)^{2}}{4}\lambda(E'_{(t-r-1)(t-r)},\vec{x}).$$
So
\begin{eqnarray}\label{eq11}
\lambda(G_{2},\vec{y})-\lambda(G_{1},\vec{x})&=& \frac{(a-b)^{2}}{4}\lambda(E'_{(t-r-1)(t-r)},\vec{x})+\lambda(G_{2},\vec{x})-\lambda(G_{1},\vec{x}) \nonumber\\
&=&\frac{(a-b)^{2}}{4}\lambda(E'_{(t-r-1)(t-r)},\vec{x})-(a-b)x_{t-r+2}x_{t-r+3}\ldots x_{t-1}x_{t}\nonumber \\
&=&\frac{(a-b)}{4}[(a-b)\lambda(E'_{(t-r-1)(t-r)},\vec{x})-4x_{t-r+2}x_{t-r+3}\ldots x_{t-1}x_{t}].
\end{eqnarray}
 By Remark \ref{r1}(b), we have
\begin{eqnarray}\label{eq12}
a-b={\lambda(E_{(t-r-1)\setminus (t-r)}) \over \lambda(E_{(t-r-1)(t-r)}) }=\frac{(x_{t-r+1}+x_{t-r+2})x_{t-r+3}\ldots x_{t-1}x_{t}}{\lambda(E_{(t-r-1)(t-r)},\vec{x})}=\frac{2x_{t-r+2}x_{t-r+3}\ldots x_{t-1}x_{t}}{\lambda(E'_{(t-r-1)(t-r)},\vec{x})}
\end{eqnarray}
since $x_{t-r+1}=x_{t-r+2}$ and $\lambda(E_{(t-r-1)(t-r)},\vec{x})=\lambda(E'_{(t-r-1)(t-r)},\vec{x})$. So
\begin{eqnarray}\label{eq13}
\lambda(G_{2},\vec{y})-\lambda(G_{1},\vec{x})=-\frac{x_{t-r+2}x_{t-r+3}\ldots x_{t-1}x_{t}}{\lambda(E'_{(t-r-1)(t-r)},\vec{x})}x_{t-r+2}x_{t-r+3}\ldots x_{t-1}x_{t}.
\end{eqnarray}
Consider a new weighting for $G_{2}$: ${\vec z}=(z_1, z_2, \ldots, z_t)$ given by $z_i=x_i$ for $i\neq t-r+1 , i\neq t-r+2$ and $z_{t-r+1}=y_{t-r+1}+\eta$, $z_{t-r+2}=y_{t-r+2}-\eta$. Then
\begin{eqnarray}
\lambda(G_{2},\vec{z})-\lambda(G_{2},\vec{y})&=& \eta[\lambda(E'_{t-r+1},\vec{y})-\lambda(E'_{t-r+2},\vec{y})]-\eta^{2}\lambda(E'_{(t-r+1)(t-r+2)},\vec{y}) \nonumber\\
&=&\eta[y_{t-r}y_{t-r+3}\ldots y_{t-1}y_{t}+y_{t-r-1}y_{t-r+3}\ldots y_{t-1}y_{t})\nonumber \\
&& -(y_{t-r+1}-y_{t-r+2})\lambda(E'_{(t-r+1)(t-r+2)},\vec{y})]-\eta^{2}\lambda(E'_{(t-r+1)(t-r+2)},\vec{y}).
\end{eqnarray}
Note that $y_{t-r+1}=y_{t-r+2}$, we have
\begin{eqnarray}\label{eq14}
\lambda(G_{2},\vec{z})-\lambda(G_{2},\vec{y})&=&\eta(y_{t-r}y_{t-r+3}\ldots
y_{t-1}y_{t}+y_{t-r-1}y_{t-r+3}\ldots y_{t-1}y_{t}) - \nonumber \\
&& - \eta^{2}\lambda(E'_{(t-r+1)(t-r+2)},\vec{y}).
\end{eqnarray}
Let
$$\eta=\frac{(y_{t-r-1}+y_{t-r})y_{t-r+3}\ldots y_{t-1}y_{t}}{2\lambda(E'_{(t-r+1)(t-r+2)},\vec{y})}.$$
Since
$$12 \ldots (r-2)\in E'_{(t-r+1)(t-r+2)},$$
then $$\eta=\frac{(y_{t-r-1}+y_{t-r})y_{t-r+3}\ldots y_{t-1}y_{t}}{2\lambda(E'_{(t-r+1)(t-r+2)},\vec{y})}\leq \frac{(y_{t-r-1}+y_{t-r})y_{t-r+3}\ldots y_{t-1}y_{t}}{2y_{1}y_{2}\ldots y_{r-2}}\leq y_{t}.$$
Hence ${\vec z}=(z_1, z_2, \ldots, z_t)$ is also a legal weighting for $G_{2}$, and
\begin{eqnarray}\label{eq15}
\lambda(G_{2},\vec{z})-\lambda(G_{2},\vec{y})&=& \frac{(y_{t-r-1}+y_{t-r})^{2}y_{t-r+3}^{2}\ldots y_{t-1}^{2}y_{t}^{2}}{2\lambda(E'_{(t-r+1)(t-r+2)},\vec{y})} - \nonumber\\
& & - [\frac{(y_{t-r-1}+y_{t-r})y_{t-r+3}\ldots y_{t-1}y_{t}}{2\lambda(E'_{(t-r+1)(t-r+2)},\vec{y})}]^{2}\lambda(E'_{(t-r+1)(t-r+2)},\vec{y}) \nonumber\\
&=& \frac{(y_{t-r-1}+y_{t-r})^{2}y_{t-r+3}^{2}\ldots y_{t-1}^{2}y_{t}^{2}}{4\lambda(E'_{(t-r+1)(t-r+2)},\vec{y})}.
\end{eqnarray}
Using (\ref{eq13}) and (\ref{eq15}), we have
\begin{eqnarray}\label{eq17}
\lambda(G_{2},\vec{z})-\lambda(G_{1},\vec{x})&=& \frac{(y_{t-r-1}+y_{t-r})^{2}y_{t-r+3}^{2}\ldots y_{t-1}^{2}y_{t}^{2}}{4\lambda(E'_{(t-r+1)(t-r+2)},\vec{y})} \nonumber\\
& & - \frac{x_{t-r+2}^{2}x_{t-r+3}^{2}\ldots x_{t-1}^{2}x_{t}^{2}}{\lambda(E'_{(t-r-1)(t-r)},\vec{x})}.
\end{eqnarray}
Note that $y_{t-r-1}+y_{t-r}=x_{t-r-1}+x_{t-r}=a+b, y_{t-r+3}=x_{t-r+3},\ldots ,y_{t}=x_{t}=c$; and $$\lambda(E'_{(t-r-1)(t-r)},\vec{x})=\lambda(E'_{(t-r-1)(t-r)},\vec{y}).$$
So
\begin{eqnarray*}
\lambda(G_{2},\vec{z})-\lambda(G_{1},\vec{x})&=& \frac{(a+b)^{2}c^{2(r-2)}}{4\lambda(E'_{(t-r+1)(t-r+2)},\vec{y})}-\frac{b^{2}c^{2(r-2)}}{\lambda(E'_{(t-r-1)(t-r)},\vec{y})}\nonumber\\
&\geq& \frac{b^{2}c^{2(r-2)}}{\lambda(E'_{(t-r+1)(t-r+2)},\vec{y})}-\frac{b^{2}c^{2(r-2)}}{\lambda(E'_{(t-r-1)(t-r)},\vec{y})}.
\end{eqnarray*}
Since $G_{2}$ is left-compressed,
we have $$\lambda(E'_{(t-r+1)(t-r+2)},\vec{y})\leq\lambda(E'_{(t-r-1)(t-r)},\vec{y}).$$
Hence $$\lambda(G_{2})\geq \lambda(G_{2},\vec{z})\geq\lambda(G_{1},\vec{x})=\lambda(G_{1}).$$
This proves the theorem.\qed

\bigskip
{\em Proof of Corollary \ref{corollary 1}.}
Let $m$ and $t$ be integers satisfying ${t-1 \choose 3} \le m \le {t \choose 3}-1.$ Let $G=(V,E)$ be a $3$-graph with $m$ edges such that $\lambda(G)=\lambda_{m}^{3}$. Applying Lemma \ref{LemmaTal8}, we can assume that $G$ is left-compressed. Let $\vec{x}=(x_{1},x_{2},\ldots ,x_{n})$  be an optimal weighting for $G$ satisfying $x_1 \ge x_2 \ge \ldots \ge x_k >x_{k+1}=\ldots=x_{n}=0$.

In \cite{T}, the following result is proved.
\begin{lemma} (Talbot \cite{T}) \label{LemmaTal9}
$$|E|\geq {k-1 \choose 3}+{k-2 \choose 2}-(k-2).$$
\end{lemma}
 We claim that $k\leq t$. Otherwise $k\geq t+1$ and Lemma \ref{LemmaTal9} implies that
\begin{eqnarray*}
m=|E|&\geq& {k-1 \choose 3}+{k-2 \choose 2}-(k-2)\nonumber\\
&\geq&{t \choose 3}+{t-1 \choose 2}-(t-1)\nonumber\\
&\ge &{t \choose 3}
\end{eqnarray*}
which contradicts to the assumption that $m={t \choose 3}-3.$
Hence $k\le t$. Combining with Theorem \ref{theorem 2}, we see that the corollary follows. \qed

\section{Proof of Theorem \ref{theorem 3} and Corollary \ref{corollary 2}}\label{proof2}

{\em Proof of Theorem \ref{theorem 3}.}  By Lemma \ref{leftcom},
we only need to consider left-compressed $r$-graphs on $[t]$  with $m={t \choose r}-4$ edges. Every left-compressed $r$-graph on $[t]$  with $m={t \choose r}-4$ edges can be obtained by removing four $r$-tuples from $[t]^{(r)}$ such that if an $r$-tuple is removed then all its ancestors should be removed. In view of Figure 1, these four $r$-tuples to be removed are
$$\{(t-r+1)(t-r+2)\ldots t,(t-r)(t-r+2)\ldots t,(t-r-1)(t-r+2)\ldots t,(t-r-2)(t-r+2)\ldots t\}$$
or
$$\{(t-r+1)(t-r+2)\ldots t,(t-r)(t-r+2)\ldots t, (t-r-1)(t-r+2)\ldots t,(t-r)(t-r+1)(t-r+3)\ldots t\}$$
or
$$\{(t-r+1)(t-r+2)\ldots t,(t-r)(t-r+2)\ldots t, (t-r)(t-r+1)(t-r+3)\ldots t,(t-r)(t-r+1)(t-r+2)(t-r+4)\ldots t\}.$$

Therefore, there are only three different left-compressed $r$-graphs with $m={t \choose r}-4$ edges on $[t]$. They are
$$G_{1}=([t],E) \ {\rm with \ the \ edge \ set }$$
$$E=[t]^{(r)}\setminus\{(t-r+1)(t-r+2)\ldots t,(t-r)(t-r+2)\ldots t,(t-r-1)(t-r+2)\ldots t,(t-r-2)(t-r+2)\ldots t\},$$
$$G_{2}=([t],E') \ {\rm with \ the \ edge \ set }$$
$$E'=[t]^{(r)}\setminus\{(t-r+1)(t-r+2)\ldots t,(t-r)(t-r+2)\ldots t, (t-r-1)(t-r+2)\ldots t,(t-r)(t-r+1)(t-r+3)\ldots t\},$$
and
$$G_{3}=([t],E'') \ {\rm with \ the \ edge \ set }$$
$$E''=[t]^{(r)}\setminus\{(t-r+1)(t-r+2)\ldots t,(t-r)(t-r+2)\ldots t, (t-r)(t-r+1)(t-r+3)\ldots t,$$
$$(t-r)(t-r+1)(t-r+2)(t-r+4)\ldots t\}.$$
Clearly, $G_{1}$ is formed by taking the first $m$ sets in the colex ordering of ${\mathbb N}^{(r)}$. So in order to prove Theorem \ref{theorem 3}, we only need to prove $\lambda(G_{1})\geq\lambda(G_{2})$ and $\lambda(G_{1})\geq\lambda(G_{3})$.

First, we show that $\lambda(G_{2})\leq \lambda(G_1)$.

Let $\vec{x}=(x_{1},x_{2},\ldots ,x_{t})$  be an optimal weighting for $G_{2}$ satisfying $x_1 \ge x_2 \ge \ldots \ge x_t \ge 0$. Note that $x_t>0$.  If $x_t=0$, then $\lambda(G_2)=\lambda([t-1]^{(r)})$. However, if we take
a legal weighting $\vec{x}=(x_1,\ldots,x_t)$, where $x_1=x_2=\cdots=x_{t-2}={1 \over t-1}$ and $x_{t-1}=x_{t}={1 \over 2(t-1)}$, then $\lambda(G_2,\vec{x} )> \lambda([t-1]^{(r)})$. This contradiction implies that $x_t>0$.
Since $G_{2}$ is left-compressed and $E'_{i\setminus j}=\emptyset$ for $i, j$ satisfying $1\le i<j\le t-r-2$, or $t-r\le i<j\le t-r+1$ or $t-r+3\le i<j\le t$,  by Remark \ref{r1}(b), we have $x_1=x_2=\cdots=x_{t-r-2}=a$, $x_{t-r-1}=b,x_{t-r}=x_{t-r+1}=c$, and $x_{t-r+2}=d$, $x_{t-r+3}=x_{t-r+4}=\ldots =x_{t}=e$. Note that
\begin{eqnarray}\label{g12xe}
\lambda(G_{1},\vec{x})-\lambda(G_{2},\vec{x})=c^2e^{r-2}-ade^{r-2}.
\end{eqnarray}
Consider a weighting for $G_{1}$: ${\vec y}=(y_1, y_2, \ldots, y_t)$ given by $y_i=x_i$ for $i\neq t-r+1$, $i\neq t-r+2$ and $y_{t-r+1}=x_{t-r+1}+\delta$, $y_{t-r+2}=x_{t-r+2}-\delta$. Then
\begin{eqnarray}\label{eq10e}
\lambda(G_{1},\vec{y})-\lambda(G_{1},\vec{x})&=& \delta[\lambda(E_{t-r+1},\vec{x})-\lambda(E_{t-r+2},\vec{x})]-\delta^{2}\lambda(E_{(t-r+1)(t-r+2)},\vec{x})\nonumber\\
&=&\delta[x_{t-r}x_{t-r+3}\ldots x_{t}+x_{t-r-1}x_{t-r+3}\ldots x_{t}+x_{t-r-2}x_{t-r+3}\ldots x_{t} \nonumber\\
&& -(x_{t-r+1}-x_{t-r+2})\lambda(E_{(t-r+1)(t-r+2)},\vec{x})]-\delta^{2}\lambda(E_{(t-r+1)(t-r+2)},\vec{x})\nonumber\\
&=&\delta[(a+b+c)e^{r-2}-(c-d)\lambda(E_{(t-r+1)(t-r+2)},\vec{x})] \nonumber\\
&& - \delta^{2}\lambda(E_{(t-r+1)(t-r+2)},\vec{x}).
\end{eqnarray}
Let
$$\delta=\frac{(a+b+c)e^{r-2}}{2\lambda(E_{(t-r+1)(t-r+2)},\vec{x})}-\frac{c-d}{2}=\frac{(a+b+c)e^{r-2}}{2\lambda(E'_{(t-r+1)(t-r+2)},\vec{x})}-\frac{c-d}{2}.$$

 By Remark \ref{r1}(b), we have
$$c-d={\lambda(E'_{(t-r+1)\setminus (t-r+2)},\vec{x}) \over \lambda(E'_{(t-r+1)(t-r+2)},\vec{x})}=\frac{be^{r-2}}{\lambda(E'_{(t-r+1)(t-r+2)},\vec{x})}.$$
So
\begin{eqnarray}\label{delta}
\delta=\frac{(a+c)e^{r-2}}{2\lambda(E'_{(t-r+1)(t-r+2)},\vec{x})}.
\end{eqnarray}
Clearly $\delta>0$. Since
$$1(t-r+4)\ldots t\in E'_{(t-r+1)(t-r+2)},$$
$$1(t-r+3)(t-r+5)\ldots t\in E'_{(t-r+1)(t-r+2)}$$
and in view of (\ref{delta}),  $\delta<e$.  So ${\vec y}=(y_1, y_2, \ldots, y_t)$ is also a legal weighting for $G_{1}$ and
\begin{eqnarray}\label{eq11e}
\lambda(G_{1},\vec{y})-\lambda(G_{1},\vec{x})=\frac{(a+c)^{2}e^{2(r-2)}}{4\lambda(E'_{(t-r+1)(t-r+2)},\vec{x})}
\end{eqnarray}
since $\lambda(E'_{(t-r+1)(t-r+2)},\vec{x})=\lambda(E_{(t-r+1)(t-r+2)},\vec{x}).$

Consider a new weighting for $G_{1}$: ${\vec z}=(z_1, z_2, \ldots, z_t)$ given by $z_i=y_i$ for $i\neq t-r-2$, $i\neq t-r$ and $z_{t-r-2}=y_{t-r-2}-\eta$, $z_{t-r}=y_{t-r}+\eta$. Then
\begin{eqnarray}\label{eq12e}
\lambda(G_{1},\vec{z})-\lambda(G_{1},\vec{y})&=& \eta[\lambda(E_{t-r},\vec{y})-\lambda(E_{t-r-2},\vec{y})]-\eta^{2}\lambda(E_{(t-r-2)(t-r)},\vec{y})\nonumber\\
&=&\eta(y_{t-r-2}-y_{t-r})\lambda(E_{(t-r-2)(t-r)},\vec{y})-\eta^{2}\lambda(E_{(t-r-2)(t-r)},\vec{y})\nonumber\\
&=&\eta(a-c)\lambda(E_{(t-r-2)(t-r)},\vec{y})-\eta^{2}\lambda(E_{(t-r-2)(t-r)},\vec{y}).
\end{eqnarray}
Let $\eta=\frac{a-c}{2}$. Clearly, ${\vec z}=(z_1, z_2, \ldots, z_t)$ is also a legal weighting for $G_{1}$, and
\begin{eqnarray}\label{eq14e}
\lambda(G_{1},\vec{z})-\lambda(G_{1},\vec{y})=\frac{(a-c)^{2}}{4}\lambda(E_{(t-r-2)(t-r)},\vec{y}).
\end{eqnarray}
Adding (\ref{g12xe}),  (\ref{eq11e}) and (\ref{eq14e}), we have
\begin{eqnarray}\label{eq15e}
\lambda(G_{1},\vec{z})-\lambda(G_{2},\vec{x})&=&\frac{(a+c)^{2}e^{2(r-2)}}{4\lambda(E'_{(t-r+1)(t-r+2)},\vec{x})}+\frac{(a-c)^{2}}{4}\lambda(E_{(t-r-2)(t-r)},\vec{y})+c^{2}e^{r-2}-ade^{r-2}\nonumber\\
&\geq&\frac{(a+c)^{2}e^{2(r-2)}}{4\lambda(E'_{(t-r+1)(t-r+2)},\vec{x})}+\frac{(a-c)^{2}}{4}\lambda(E_{(t-r-2)(t-r)},\vec{y})-(a-c)de^{r-2}.
\end{eqnarray}
By Remark \ref{r1}(b), we have
$$a-c={\lambda(E'_{(t-r-2)\setminus (t-r)},\vec{x}) \over \lambda(E'_{(t-r-2)(t-r)},\vec{x})}=\frac{(d+c)e^{r-2}}{\lambda(E'_{(t-r-2)(t-r)},\vec{x})}.$$
Hence
\begin{eqnarray}\label{eq16e}
\lambda(G_{1},\vec{z})-\lambda(G_{2},\vec{x})&\ge&\frac{(a+c)^{2}e^{2(r-2)}}{4\lambda(E'_{(t-r+1)(t-r+2)},\vec{x})}+\frac{(d+c)^{2}e^{2(r-2)}}{4(\lambda(E'_{(t-r-2)(t-r)},\vec{x}))^{2}}\lambda(E_{(t-r-2)(t-r)},\vec{y}) \nonumber\\
&& - \frac{(d+c)de^{2(r-2)}}{\lambda(E'_{(t-r-2)(t-r)},\vec{x})}\nonumber\\
&\geq&\frac{(a+c)^{2}e^{2(r-2)}}{4\lambda(E'_{(t-r+1)(t-r+2)},\vec{x})}+\frac{(d+c)^{2}e^{2(r-2)}}{4(\lambda(E'_{(t-r-2)(t-r)},\vec{x}))^{2}}\lambda(E_{(t-r+1)(t-r+2)},\vec{y}) \nonumber\\
&& -\frac{(d+c)^2e^{2(r-2)}}{2\lambda(E'_{(t-r-2)(t-r)},\vec{x})}\nonumber\\
&=&\frac{(a+c)^{2}e^{2(r-2)}}{4\lambda(E'_{(t-r+1)(t-r+2)},\vec{x})}+\frac{(d+c)^{2}e^{2(r-2)}}{4(\lambda(E'_{(t-r-2)(t-r)},\vec{x}))^{2}}\lambda(E_{(t-r+1)(t-r+2)},\vec{x}) \nonumber\\
&& -\frac{(d+c)^2e^{2(r-2)}}{2\lambda(E'_{(t-r-2)(t-r)},\vec{x})}\nonumber\\
&=&\frac{(a+c)^{2}e^{2(r-2)}}{4\lambda(E'_{(t-r+1)(t-r+2)},\vec{x})}+\frac{(d+c)^{2}e^{2(r-2)}}{4(\lambda(E'_{(t-r-2)(t-r)},\vec{x}))^{2}}\lambda(E'_{(t-r+1)(t-r+2)},\vec{x})  \nonumber\\
&& - \frac{(d+c)^2e^{2(r-2)}}{2\lambda(E'_{(t-r-2)(t-r)},\vec{x})}\nonumber\\
\end{eqnarray}
Observe  that
\begin{eqnarray}\label{eq166e}
&&[\frac{(d+c)^2e^{2(r-2)}\lambda(E'_{(t-r+1)(t-r+2)},\vec{x})}{4(\lambda(E'_{(t-r-2)(t-r)},\vec{x}))^{2}}
-\frac{(d+c)^2e^{2(r-2)}}{2\lambda(E'_{(t-r-2)(t-r)},\vec{x})}] \nonumber \\
&&-[\frac{(d+c)^2e^{2(r-2)}\lambda(E'_{(t-r+1)(t-r+2)},\vec{x})}{4(\lambda(E'_{(t-r+1)(t-r+2)},\vec{x}))^{2}}-\frac{(d+c)^2e^{2(r-2)}}{2\lambda(E'_{(t-r+1)(t-r+2)},\vec{x})}] \nonumber \\
&=&(d+c)^2e^{2(r-2)}[\lambda(E'_{(t-r+1)(t-r+2)},\vec{x})({1 \over 4(\lambda(E'_{(t-r-2)(t-r)},\vec{x}))^{2}}-{1 \over 4(\lambda(E'_{(t-r+1)(t-r+2)},\vec{x}))^{2}}) \nonumber \\
&&-(\frac{1}{2\lambda(E'_{(t-r-2)(t-r)},\vec{x})}-\frac{1}{2\lambda(E'_{(t-r+1)(t-r+2)},\vec{x})})]\nonumber \\
&=&(d+c)^2e^{2(r-2)}(\frac{1}{2\lambda(E'_{(t-r-2)(t-r)},\vec{x})}-\frac{1}{2\lambda(E'_{(t-r+1)(t-r+2)},\vec{x})}) \nonumber \\
&& \times [\lambda(E'_{(t-r+1)(t-r+2)},\vec{x})(\frac{1}{2\lambda(E'_{(t-r-2)(t-r)},\vec{x})}+\frac{1}{2\lambda(E'_{(t-r+1)(t-r+2)},\vec{x})})-1]\nonumber \\
&\ge& 0 .
\end{eqnarray}
The last inequality is true because of the following:  since $G$ is left-compressed, then $\lambda(E'_{(t-r+1)(t-r+2)},\vec{x})\le \lambda(E'_{(t-r-2)(t-r)},\vec{x})$. So  $\frac{1}{2\lambda(E'_{(t-r-2)(t-r)},\vec{x})}-\frac{1}{2\lambda(E'_{(t-r+1)(t-r+2)},\vec{x})}\le 0$ and $\lambda(E'_{(t-r+1)(t-r+2)},\vec{x})(\frac{1}{2\lambda(E'_{(t-r-2)(t-r)},\vec{x})}+\frac{1}{2\lambda(E'_{(t-r+1)(t-r+2)},\vec{x})})-1\le 0$.

Combining (\ref{eq16e}) and (\ref{eq166e}), we have
\begin{eqnarray*}
\lambda(G_{1},\vec{z})-\lambda(G_{2},\vec{x})&\geq&\frac{(a+c)^{2}e^{2(r-2)}}{4\lambda(E'_{(t-r+1)(t-r+2)},\vec{x})}+\frac{(d+c)^{2}e^{2(r-2)}}{4(\lambda(E'_{(t-r+1)(t-r+2)},\vec{x}))^{2}}\lambda(E'_{(t-r+1)(t-r+2)},\vec{x})  \nonumber\\
&& - \frac{(d+c)^2e^{2(r-2)}}{2\lambda(E'_{(t-r+1)(t-r+2)},\vec{x})}\nonumber\\
&\geq& 0.
\end{eqnarray*}
Hence
$$\lambda(G_{1})\geq\lambda(G_{1},\vec{z})\geq\lambda(G_{2},\vec{x})=\lambda(G_{2}).$$

Next we prove $\lambda(G_{3})\leq \lambda(G_1)$.

Let $\vec{x'}=(x'_{1},x'_{2},\ldots ,x'_{t})$ be an optimal weighting for $G_{3}$ satisfying $x'_1 \ge x'_2 \ge \ldots \ge x'_t \ge 0$. Note that $x'_t>0$.  If $x'_t=0$, then $\lambda(G_3)=\lambda([t-1]^{(r)})$. However, if we take
a legal weighting $\vec{x'}=(x'_1,\ldots,x'_t)$, where $x'_1=x'_2=\cdots=x'_{t-2}={1 \over t-1}$ and $x'_{t-1}=x'_{t}={1 \over 2(t-1)}$, then $\lambda(G_3,\vec{x'} )> \lambda([t-1]^{(r)})$. This contradiction implies that $x'_t>0$.
Since $G_{3}$ is left-compressed and $E''_{i\setminus j}=\emptyset$ for $i, j$ satisfying $1\le i<j\le t-r-1$, or $t-r\le i<j\le t-r+3$ or $t-r+4\le i<j\le t$, by Remark \ref{r1}(b), we have $x'_1=x'_2=\cdots=x'_{t-r-1}=a'$, $x'_{t-r}=x'_{t-r+1}=x'_{t-r+2}=x'_{t-r+3}=b'$, $x'_{t-r+4}=x'_{t-r+5}=\cdots =x'_{t}=c'$. Note that
\begin{eqnarray}\label{g13x}
\lambda(G_{1},\vec{x'})-\lambda(G_{3},\vec{x'})=(2b'^3-2a'b'^2)c'^{r-3}.
\end{eqnarray}
By Remark \ref{r1}(b), we have
$$a'-b'={\lambda(E''_{(t-r-1)\setminus (t-r)},\vec{x}) \over \lambda(E''_{(t-r-1)(t-r)},\vec{x})}=\frac{3b'^2c'^{t-3}}{\lambda(E''_{(t-r-1)(t-r)},\vec{x})}\le c'\le b'.$$ Therefore,
\begin{equation}\label{ale2b}
a'\le 2b'.
\end{equation}
Consider a weighting for $G_{1}$: $\vec {y'}=(y'_1, y'_2, \ldots, y'_t)$ given by $y'_i=x'_i$ for $i\neq t-r$, $i\neq t-r+3$ and $y'_{t-r}=x'_{t-r}+\frac{a'-b'}{3}=\frac{a'+2b'}{3}$, $y'_{t-r+3}=x'_{t-r+3}-\frac{a'-b'}{3}=\frac{4b'-a'}{3}$. Clearly, $\vec{ y'}=(y'_1, y'_2, \ldots, y'_t)$ is a legal weighting for $G_{1}$, and
\begin{eqnarray}\label{eq18e}
\lambda(G_{1},\vec{y'})-\lambda(G_{1},\vec{x'})&=&\frac{a'-b'}{3}[\lambda(E_{t-r},\vec{x'})-\lambda(E_{t-r+3},\vec{x'})]-(\frac{a'-b'}{3})^{2}\lambda(E_{(t-r)(t-r+3)},\vec{x'})\nonumber\\
&=&\frac{a'-b'}{3}[x'_{t-r+1}x'_{t-r+2}x'_{t-r+4}\ldots x'_{t}+x'_{t-r-1}x'_{t-r+2}x'_{t-r+4}\ldots x'_{t} + \nonumber\\
&& + x'_{t-r-2}x'_{t-r+2}x'_{t-r+4}\ldots x'_{t-r}+(x'_{t-r+3}-x'_{t-r})\lambda(E_{(t-r)(t-r+3)},\vec{x'})] \nonumber\\
&& - (\frac{a'-b'}{3})^{2}\lambda(E_{(t-r)(t-r+3},\vec{x'})\nonumber\\
&=&\frac{a'-b'}{3}(b'^2c^{r-3}+2a'b'c^{r-3})-\frac{(a'-b')^2}{9}\lambda(E_{(t-r)(t-r+3)},\vec{x'}).
\end{eqnarray}
Consider a new weighting for $G_{1}$: $\vec{z'}=(z'_1, z'_2, \ldots, z'_t)$ given by $z'_i=y'_i$ for $i\neq t-r+1 , i\neq t-r+2$ and $z'_{t-r+1}=y'_{t-r+1}+\frac{a'-b'}{3}=\frac{a'+2b'}{3}$, $z'_{t-r+2}=y'_{t-r+2}-\frac{a'-b'}{3}=\frac{4b'-a'}{3}$.  Clearly $\vec{ z'}=(z'_1, z'_2, \ldots, z'_t)$ is also a legal weighting for $G_{1}$, and
\begin{eqnarray}\label{eq19e}
\lambda(G_{1},\vec{z'})-\lambda(G_{1},\vec{y'})&=& \frac{a'-b'}{3}[\lambda(E_{t-r+1},\vec{y'})-\lambda(E_{t-r+2)},\vec{y'})]-(\frac{a'-b'}{3})^{2}\lambda(E_{(t-r+1)(t-r+2)},\vec{y'})\nonumber\\
&=&\frac{a'-b'}{3}[y'_{t-r}y'_{t-r+3}\cdots y'_{t}+y'_{t-r-1}y'_{t-r+3}\cdots y'_{t}+y'_{t-r-2}y'_{t-r+3}\cdots y'_{t}  \nonumber\\
&& - (y'_{t-r+1}-y'_{t-r+2})\lambda(E_{(t-r+1)(t-r+2)},\vec{y'})]-(\frac{a'-b'}{3})^{2}\lambda(E_{(t-r+1)(t-r+2)},\vec{y'})\nonumber\\
&=&\frac{a'-b'}{3}c'^{r-3}(\frac{a'+2b'}{3}\frac{4b'-a'}{3}+a'\frac{4b'-a'}{3}+a'\frac{4b'-a'}{3})  \nonumber\\
&& - \frac{(a'-b')^2}{9}\lambda(E_{(t-r+1)(t-r+2)},\vec{y'}).
\end{eqnarray}
Again consider a new weighting for $G_{1}$: $\vec{ y''}=(y''_1, y''_2, \ldots, y''_t)$ given by $y''_i=z'_i$ for $i\neq t-r-1$, $i\neq t-r$ and $y''_{t-r-1}=z'_{t-r-1}-\frac{a'-b'}{3}=\frac{2a'+b'}{3}$, $y''_{t-r}=z'_{t-r}+\frac{a'-b'}{3}=\frac{2a'+b'}{3}$. Clearly, $\vec{ y''}=(y''_1, y''_2, \ldots, y''_t)$ is also a legal weighting for $G_{1}$, and
\begin{eqnarray}\label{eq20e}
\lambda(G_{1},\vec{y''})-\lambda(G_{1},\vec{z'})&=&\frac{a'-b'}{3}[\lambda(E_{t-r},\vec{z'})-\lambda(E_{t-r-1},\vec{z'})]-(\frac{a'-b'}{3})^{2}\lambda(E_{(t-r-1)(t-r)},\vec{z'})\nonumber\\
&=&\frac{a'-b'}{3}(a'-\frac{a'+2b'}{3})\lambda(E_{(t-r-1)(t-r)},\vec{z'})-(\frac{a'-b'}{3})^{2}\lambda(E_{(t-r-1)(t-r)},\vec{z'})\nonumber\\
&=&\frac{(a'-b')^2}{9}\lambda(E_{(t-r-1)(t-r)},\vec{z'}).
\end{eqnarray}
Consider a new weighting for $G_{1}$ once more: $\vec{ z''}=(z''_1, z''_2, \ldots, z''_t)$ given by $z''_i=y''_i$ for $i\neq t-r-2$, $i\neq t-r+1$ and $z''_{t-r-2}=y''_{t-r-2}-\frac{a'-b'}{3}=\frac{2a'+b'}{3}$, $z''_{t-r+1}=y''_{t-r+1}+\frac{a'-b'}{3}=\frac{2a'+b'}{3}$. Clearly, $\vec{ z''}=(z''_1, z''_2, \ldots, z''_t)$ is also a legal weighting for $G_{1}$, and
\begin{eqnarray}\label{eq21e}
\lambda(G_{1},\vec{z''})-\lambda(G_{1},\vec{y''})&=&\frac{a'-b'}{3}[\lambda(E_{t-r+1},\vec{y''})-\lambda(E_{t-r-2},\vec{y''})]-(\frac{a'-b'}{3})^{2}\lambda(E_{(t-r-2)(t-r+1)},\vec{y''})\nonumber\\
&=&\frac{a'-b'}{3}(a'-\frac{a'+2b'}{3})\lambda(E_{(t-r-2)(t-r+1)},\vec{y''})-(\frac{a'-b'}{3})^{2}\lambda(E_{(t-r-2)(t-r+1)},\vec{y''})\nonumber\\
&=&\frac{(a'-b')^2}{9}\lambda(E_{(t-r-2)(t-r+1)},\vec{y''}).
\end{eqnarray}
Adding (\ref{g13x}), (\ref{eq18e}), (\ref{eq19e}),(\ref{eq20e}),and (\ref{eq21e}), we have
\begin{eqnarray}\label{eq22e}
\lambda(G_{1},\vec{z''})-\lambda(G_{3},\vec{x'})&=&\frac{(a'-b')^{2}}{9}[\lambda(E_{(t-r-1)(t-r)},\vec{z'})+\lambda(E_{(t-r-2)(t-r+1)},\vec{y''}) - \nonumber\\
&&  \lambda(E_{(t-r)(t-r+3)},\vec{x'})-\lambda(E_{(t-r+1)(t-r+2)},\vec{y'})]+\frac{a'-b'}{3}c'^{r-3}(b'^2+2a'b') \nonumber\\
&& - \frac{a'-b'}{3}c'^{r-3}(\frac{a'+2b'}{3}\frac{4b'-a'}{3}+a'\frac{4b'-a'}{3}+a'\frac{4b'-a'}{3}) \nonumber \\
&& + (2b^3-2a'b'^2)c'^{r-3}. 
\end{eqnarray}
Note that
\begin{eqnarray}\label{eq23e}
&&\lambda(E_{(t-r-2)(t-r+1)},\vec{y''})-\lambda(E_{(t-r-2)(t-r+1)},\vec{z'})\nonumber\\
&=&\frac{a'-b'}{3}[\lambda(E_{(t-r-2)(t-r)(t-r+1)},\vec{z'})
-\lambda(E_{(t-r-2)(t-r-1)(t-r+1)},\vec{z'})] \nonumber\\
&& - \frac{(a'-b')^2}{9}\lambda(E_{(t-r-2)(t-r-1)(t-r)(t-r+1)},\vec{z'})\nonumber\\
&=&\frac{a'-b'}{3}(z'_{t-r-1}-z'_{t-r})\lambda(E_{(t-r-2)(t-r-1)(t-r)(t-r+1)},\vec{z'}) \nonumber\\
&& -\frac{(a'-b')^2}{9}\lambda(E_{(t-r-2)(t-r-1)(t-r)(t-r+1)},\vec{z'})\nonumber\\
&=&\frac{(a'-b')^2}{9}\lambda(E_{(t-r-2)(t-r-1)(t-r)(t-r+1)},\vec{z'})\geq 0,
\end{eqnarray}
also note that
\begin{eqnarray}\label{eq24e}
\lambda(E_{(t-r-2)(t-r+1)},\vec{z'})\geq \lambda(E_{(t-r-2)(t-r+2)},\vec{z'})
\end{eqnarray}
since $G_1$ is left-compressed;
and note that
\begin{eqnarray}\label{eq25e}
\lambda(E_{(t-r-2)(t-r+2)},\vec{z'})-\lambda(E_{(t-r-2)(t-r+2)},\vec{y'})=\frac{a'-b'}{3}\lambda(E_{(t-r-2)(t-r+1)(t-r+2)},\vec{y'})\geq 0.
\end{eqnarray}
Using (\ref{eq23e}), (\ref{eq24e}), and (\ref{eq25e}), we have
$$\lambda(E_{(t-r-2)(t-r+1)},\vec{y''})\ge \lambda(E_{(t-r-2)(t-r+2)},\vec{y'}).$$
Therefore,
\begin{eqnarray}\label{eq26e}
&&\lambda(E_{(t-r-1)(t-r)},\vec{z'})+\lambda(E_{(t-r-2)(t-r+1)},\vec{y''})  \nonumber\\
&& - \lambda(E_{(t-r)(t-r+3)},\vec{x'})-\lambda(E_{(t-r+1)(t-r+2)},\vec{y'})\nonumber\\
&\geq &\lambda(E_{(t-r+1)(t-r+2)},\vec{z'})+\lambda(E_{(t-r-2)(t-r+2)},\vec{y'}) \nonumber\\
&& - \lambda(E_{(t-r)(t-r+3)},\vec{y'})-\lambda(E_{(t-r+1)(t-r+2)},\vec{y'}) \ge 0
\end{eqnarray}
 since $\lambda(E_{(t-r+1)(t-r+2)},\vec{z'})=\lambda(E_{(t-r+1)(t-r+2)},\vec{y'})$.

Recall that  $a'\leq 2b'$.  So
\begin{eqnarray}\label{eq27e}
&&\frac{a'-b'}{3}c'^{r-3}(b'^2+2a'b')+\frac{a'-b'}{3}c'^{r-3}(\frac{a'+2b'}{3}\frac{4b'-a'}{3}
+a'\frac{4b'-a'}{3}+a'\frac{4b'-a'}{3})+ \nonumber \\
&& +(2b'^3-2a'b'^2)c'^{r-3}\nonumber\\
&=&\frac{a'-b'}{3}c'^{r-3}(b'^2+2a'b'+\frac{a'+2b'}{3}\frac{4b'-a'}{3}+2a'\frac{4b'-a'}{3}-6b'^{2})\nonumber\\
&=&\frac{a'-b'}{27}c'^{r-3}(44a'b'-37b'^2-7a'^2)\nonumber\\
&=&\frac{(a'-b')^2}{27}c'^{r-3}(37b'-7a')\geq 0.
\end{eqnarray}
Combing (\ref{eq22e}), (\ref{eq26e}) and  (\ref{eq27e}), we have
$$\lambda(G_{1})\geq\lambda(G_{1},\vec{z''})\geq\lambda(G_{3},\vec{x'})=\lambda(G_{3}).$$
This proves the theorem.\qed

\bigskip

{\em Proof of Corollary \ref{corollary 2}.}
Let $m$ and $t$ be integers satisfying ${t-1 \choose 3} \le m \le {t \choose 3}-1.$ Let $G=(V,E)$ be a $3$-graph with $m$ edges such that $\lambda(G)=\lambda_{m}^{3}$. Applying Lemma \ref{LemmaTal8}, we can assume that $G$ is left-compressed. Let $\vec{x}=(x_{1},x_{2},\ldots ,x_{n})$  be an optimal weighting for $G$ satisfying $x_1 \ge x_2 \ge \ldots \ge x_k >x_{k+1}=\ldots=x_{n}=0$.

 We claim that $k\leq t$. Otherwise $k\geq t+1$ and Lemma\ref{LemmaTal9} implies that
\begin{eqnarray*}
m=|E|&\geq& {k-1 \choose 3}+{k-2 \choose 2}-(k-2)\nonumber\\
&\geq&{t \choose 3}+{t-1 \choose 2}-(t-1)\nonumber\\
&\ge &{t \choose 3}
\end{eqnarray*}
which contradicts to the assumption that $m={t \choose 3}-4.$
Hence $k\le t$. Combining with Theorem \ref{theorem 3}, we see that  the corollary holds. \qed

\begin{remark} 
If  a result similar to Lemma \ref{LemmaTal9} holds for $r\ge 4$, then combining with Theorems \ref{theorem 2} and  \ref{theorem 3},  Conjecture \ref{conjecture} holds for general $r\ge 4$ when $m= {t \choose r}-3$ or $m={t \choose r}-4$.

\end{remark}

\end{document}